\documentclass{amsart}

\newtheorem{theorem}[equation]{Theorem}
\newtheorem{lemma}[equation]{Lemma}
\newtheorem{corollary}[equation]{Corollary}
\newtheorem{proposition}[equation]{Proposition}

\numberwithin{equation}{section}

\newcommand{\updot}{\textstyle\cdot}

\usepackage{amscd}

\begin{document}

\title{Exponential sums on ${\bf A}^n$}
\author{Alan Adolphson}
\address{Department of Mathematics\\
Oklahoma State University\\
Stillwater, Oklahoma 74078}
\email{adolphs@math.okstate.edu}
\thanks{The first author was supported in part by NSA Grant
  \#MDA904-97-1-0068} 
\author{Steven Sperber}
\address{School of Mathematics\\
University of Minnesota\\
Minneapolis, Minnesota 55455}
\email{sperber@math.umn.edu}
\date{April 22, 1998}
\keywords{Exponential sum, $p$-adic cohomology, $l$-adic cohomology}
\subjclass{Primary 11L07, 11T23, 14F20, 14F30}
\begin{abstract}
We discuss exponential sums on affine space from the point of view of
Dwork's $p$-adic cohomology theory.
\end{abstract}
\maketitle

\section{Introduction}

Let $p$ be a prime number, $q=p^a$, ${\bf F}_q$ the finite field of
$q$ elements.  Associated to a polynomial $f\in{\bf
  F}_q[x_1,\ldots,x_n]$ and a nontrivial additive character $\Psi:{\bf
  F}_q\rightarrow{\bf C}^{\times}$ are exponential sums
\begin{equation}
S({\bf A}^n({\bf F}_{q^i}),f)=\sum_{x_1,\ldots,x_n\in{\bf F}_{q^i}}
\Psi({\rm Trace}_{{\bf F}_{q^i}/{\bf F}_q}f(x_1,\ldots,x_n))
\end{equation}
and an $L$-function
\begin{equation} L({\bf A}^n,f;t)=\exp\biggl(\sum_{i=1}^{\infty}
  S({\bf A}^n({\bf F}_{q^i}),f)\frac{t^i}{i}\biggr).
\end{equation}
Let $d=\text{degree of $f$}$ and write
\[ f=f^{(d)}+f^{(d-1)}+\cdots+f^{(0)}, \]
where $f^{(j)}$ is homogeneous of degree $j$.  A by now classical
theorem of Deligne\cite[Th\'{e}or\`{e}me 8.4]{DE1}  says that if
$(p,d)=1$ and $f^{(d)}=0$ defines a smooth hypersurface in ${\bf
  P}^{n-1}$, then $L({\bf A}^n,f;t)^{(-1)^{n+1}}$ is a polynomial of
degree $(d-1)^n$, all of whose reciprocal roots have absolute value
equal to $q^{n/2}$.  This implies the estimate
\begin{equation} |S({\bf A}^n({\bf F}_{q^i}),f)|\leq (d-1)^nq^{ni/2}.
\end{equation}

In this article, we give a $p$-adic proof of the fact that $L({\bf
  A}^n,f;t)^{(-1)^{n+1}}$ is a polynomial of degree $(d-1)^n$
(equation (2.14) and Theorem 3.8) and give $p$-adic estimates for its
reciprocal roots, namely, we find a lower bound for the $p$-adic
Newton polygon of $L({\bf A}^n,f;t)^{(-1)^{n+1}}$ (Theorem 4.3).
Using general results of Deligne\cite{DE2}, this information can be
used to compute $l$-adic cohomology and hence again obtain the
archimedian estimate (1.3) (Theorem 5.3).

For Theorems 3.8 and 4.3, we need to assume only that $\{\partial
f^{(d)}/\partial x_i\}_{i=1}^n$ form a regular sequence in ${\bf
  F}_q[x_1,\ldots,x_n]$ (or, equivalently, that $\{\partial
f^{(d)}/\partial x_i\}_{i=1}^n$ have no common zero in ${\bf
  P}^{n-1}$).  When $(p,d)=1$, this is equivalent to 
Deligne's hypothesis.  When $d$ is divisible by $p$, there are only a
few cases satisfying this regular sequence condition.  We check them
by hand in section 6 to prove the following slight generalization of
Deligne's result.
\begin{theorem}
Suppose $\{\partial f^{(d)}/\partial x_i\}_{i=1}^n$ form a regular
  sequence in ${\bf F}_q[x_1,\ldots,x_n]$.  Then $L({\bf
  A}^n,f;t)^{(-1)^{n+1}}$ is a polynomial of degree $(d-1)^n$, all of
  whose reciprocal roots have absolute value equal to $q^{n/2}$.
\end{theorem}

In the article \cite{AS1}, we dealt with exponential sums on tori.
After a general coordinate change, one can, by using the standard
toric decomposition of ${\bf A}^n$, deduce most of the results of this
article from results in \cite{AS1}.  Our main purpose here is to
develop some new methods that will be more widely applicable.  For
instance, recent results of Garc\'{\i}a\cite{GA} on exponential sums
on ${\bf A}^n$ do not seem to follow from~\cite{AS1}.

In contrast with \cite{AS1}, we work systematically with spaces of
type $C(b)$ (convergent series on a closed disk) and avoid spaces of
type $L(b)$ (bounded series on an open disk).  This ties together more
closely the calculation of $p$-adic cohomology and the estimation of
the Newton polygon of the characteristic polynomial of Frobenius,
eliminating much of section 3 of \cite{AS1}.

Another new feature of this work is the use of the spectral sequence
associated to the filtration by $p$-divisibility on the complex
$\Omega^{\updot}_{C(b)}$ (section 3 below).  Although the behavior of
this spectral sequence is rather simple in the setting of this article
(namely, $E_1^{r,s}=E_{\infty}^{r,s}$ for all $r$ and $s$), we believe
it will play a significant role in more general situations, such as
that of Garc\'{\i}a\cite{GA}.  We hope the methods developed here will
allow us to extend the results of this article to those situations.

\section{Preliminaries}

In this section, we review the results from Dwork's $p$-adic
cohomology theory that will be used in this paper.  

Let ${\bf Q}_p$ be the field of $p$-adic numbers,
$\zeta_p$ a primitive $p$-th root of unity, and $\Omega_1={\bf
  Q}_p(\zeta_p)$.  The field $\Omega_1$ is a totally ramified
extension of ${\bf Q}_p$ of degree $p-1$.  Let $K$ be the unramified
extension of ${\bf Q}_p$ of degree $a$.  Set $\Omega_0=K(\zeta_p)$.
The Frobenius automorphism $x\mapsto x^p$ of ${\rm Gal}({\bf F}_q/{\bf
  F}_p)$ lifts to a generator $\tau$ of ${\rm
  Gal}(\Omega_0/\Omega_1)(\simeq {\rm Gal}(K/{\bf Q}_p))$ by requiring
$\tau(\zeta_p)=\zeta_p$.  Let $\Omega$ be the completion of an
algebraic closure of $\Omega_0$.  Denote by ``ord'' the additive
valuation on $\Omega$ normalized by ${\rm ord}\;p=1$ and by ``${\rm
  ord}_q$'' the additive valuation normalized by ${\rm ord}_q\;q=1$.

Let $E(t)$ be the Artin-Hasse exponential series:
\[ E(t)=\exp\biggl(\sum_{i=0}^{\infty} \frac{t^{p^i}}{p^i}\biggr). \]
Let $\gamma\in\Omega_1$ be a solution of $\sum_{i=0}^{\infty}
t^{p^i}/p^i=0$ satisfying ${\rm ord}\;\gamma=1/(p-1)$ and consider
\begin{equation}
\theta(t)=E(\gamma t)=\sum_{i=0}^{\infty}
\lambda_it^i\in\Omega_1[[t]].
\end{equation}
The series $\theta(t)$ is a splitting function in Dwork's
terminology\cite{DW1}.  Furthermore, its coefficients satisfy
\begin{equation}
{\rm ord}\;\lambda_i\geq i/(p-1).
\end{equation}

We consider the following spaces of $p$-adic functions.  Let $b$ be a
positive rational number and choose a positive integer $M$ such that
$Mb/p$ and $Md/(p-1)$ are integers.  Let $\pi$ be such that
\begin{equation}
\pi^{Md}=p
\end{equation}
and put $\tilde{\Omega}_1=\Omega_1(\pi)$,
$\tilde{\Omega}_0=\Omega_0(\pi)$.  The element $\pi$ is a uniformizing
parameter for the rings of integers of $\tilde{\Omega}_1$ and
$\tilde{\Omega}_0$.  We extend $\tau\in {\rm Gal}(\Omega_0/\Omega_1)$
to a generator of ${\rm Gal}(\tilde{\Omega}_0/\tilde{\Omega}_1)$ by
requiring $\tau(\pi)=\pi$.  For $u=(u_1,\ldots,u_n)\in{\bf R}^n$, we
put $|u|=u_1+\cdots+u_n$.  Define
\begin{equation}
C(b)=\biggl\{ \sum_{u\in{\bf N}^n} A_u\pi^{Mb|u|}x^u \mid
  \text{$A_u\in\tilde{\Omega}_0$ and $A_u\rightarrow 0$ as
  $u\rightarrow\infty$} \biggr\}. 
\end{equation}
For $\xi=\sum_{u\in{\bf N}^n}A_u\pi^{Mb|u|}x^u\in C(b)$, define
\[ {\rm ord}\;\xi=\min_{u\in{\bf N}^n} \{{\rm ord}\;A_u\}. \]
Given $c\in {\bf R}$, we put
\[ C(b,c)=\{\xi\in C(b)\mid {\rm ord}\;\xi\geq c\}. \]

Let $\hat{f}=\sum_u \hat{a}_ux^u\in K[x_1,\ldots,x_n]$ be the
Teichm\"{u}ller lifting of the polynomial $f\in{\bf
  F}_q[x_1,\ldots,x_n]$, i.~e., $(\hat{a}_u)^q=\hat{a}_u$ and 
the reduction of $\hat{f}$ modulo $p$ is $f$.  Set
\begin{align}
F(x)& = \prod_u \theta(\hat{a}_ux^u), \\
F_0(x)& = \prod_{i=0}^{a-1}\prod_u \theta((\hat{a}_ux^u)^{p^i}).
\end{align}
The estimate (2.2) implies that $F\in
C(b,0)$ for all $b<1/(p-1)$ and $F_0\in C(b,0)$ for all $b<p/q(p-1)$.
Define an operator $\psi$ on formal power series by
\begin{equation}
\psi\biggl(\sum_{u\in{\bf N}^n}A_ux^u\biggr)=\sum_{u\in{\bf N}^n}
A_{pu}x^u.
\end{equation}
It is clear that $\psi(C(b,c))\subseteq C(pb,c)$.
For $0<b<p/(p-1)$, let $\alpha=\psi^a\circ F_0$ be the composition
\[ C(b)\hookrightarrow C(b/q)\xrightarrow{F_0}
C(b/q)\xrightarrow{\psi^a} C(b). \] 
Then $\alpha$ is a completely continuous $\tilde{\Omega}_0$-linear
endomorphism of $C(b)$.  We shall also need to consider
$\beta=\tau^{-1}\circ\psi\circ F$, which is a completely continuous
$\tilde{\Omega}_1$-linear (or $\tilde{\Omega}_0$-semilinear)
endomorphism of $C(b)$.  Note that $\alpha=\beta^a$.

Set $\hat{f}_i=\partial\hat{f}/\partial x_i$ and let
$\gamma_l=\sum_{i=0}^l \gamma^{p^i}/p^i$.  By the definition of
$\gamma$, we have 
\begin{equation}
{\rm ord}\;\gamma_l\geq \frac{p^{l+1}}{p-1}-l-1.
\end{equation}
For $i=1,\ldots,n$, define differential operators $D_i$ by
\begin{equation}
D_i=\frac{\partial}{\partial x_i}+H_i,
\end{equation}
where 
\begin{equation}
H_i=\sum_{l=0}^{\infty} \gamma_lp^lx_i^{p^l-1}\hat{f}_i^{\tau^l}(x^{p^l}) \in
C\biggl(b,\frac{1}{p-1}-b\frac{d-1}{d}\biggr)
\end{equation}
for $b<p/(p-1)$.  Thus $D_i$ and ``multiplication by $H_i$'' operate
on $C(b)$ for $b<p/(p-1)$.  

To understand the definition of the $D_i$, put
\begin{align*}
\hat{\theta}(t)&=\prod_{i=0}^{\infty}\theta(t^{p^i}), \\
\hat{F}(x)&=\prod_u \hat{\theta}(\hat{a}_ux^u),
\end{align*}
so that
\begin{align*}
F(x)&=\hat{F}(x)/\hat{F}(x^p), \\
F_0(x)&=\hat{F}(x)/\hat{F}(x^q).
\end{align*}
Then formally
\begin{align*}
\alpha&=\hat{F}(x)^{-1}\circ\psi^a\circ\hat{F}(x) \\
\beta&=\hat{F}(x)^{-1}\circ\tau^{-1}\circ\psi\circ\hat{F}(x).
\end{align*}
It is trivial to check that $x_i\partial/\partial x_i$ and $\psi$
commute up to a factor of $p$, hence the differential operators
\[ \hat{F}^{-1}\circ x_i\frac{\partial}{\partial x_i}\circ\hat{F} =
x_i\frac{\partial}{\partial x_i}+\frac{x_i\partial\hat{F}/\partial
  x_i}{\hat{F}} \]
formally commute with $\alpha$ (up to a factor of $q$) and $\beta$ (up
  to a factor of $p$).  From the definitions, one gets
\[
  \hat{\theta}(t)=\exp\biggl(\sum_{l=0}^{\infty}\gamma_lt^{p^l}\biggr). \]
It then follows that
\[ \frac{x_i\partial\hat{F}/\partial x_i}{\hat{F}}=x_iH_i, \]
which gives
\begin{align}
\alpha\circ x_iD_i&=qx_iD_i\circ\alpha, \\
\beta\circ x_iD_i&=px_iD_i\circ\beta.
\end{align}

Consider the de Rham-type complex $(\Omega_{C(b)}^{\updot},D)$, where
\[ \Omega_{C(b)}^k=\bigoplus_{1\leq i_1<\cdots<i_k\leq n} C(b)\,
dx_{i_1}\wedge\cdots\wedge dx_{i_k} \]
and $D:\Omega_{C(b)}^k\rightarrow \Omega_{C(b)}^{k+1}$ is defined by
\[ D(\xi\,dx_{i_1}\wedge\cdots\wedge dx_{i_k})=\biggl(\sum_{i=1}^n
D_i(\xi)\,dx_i\biggr)\wedge dx_{i_1}\wedge\cdots\wedge dx_{i_k}. \]
We extend the mapping $\alpha$ to a mapping
$\alpha_{\updot}:\Omega_{C(b)}^{\updot}\rightarrow
\Omega_{C(b)}^{\updot}$ defined by linearity and the formula
\[ \alpha_k(\xi\,dx_{i_1}\wedge\cdots\wedge dx_{i_k})=
q^{n-k}\frac{1}{x_{i_1}\cdots x_{i_k}}\alpha(x_{i_1}\cdots
x_{i_k}\xi)\, dx_{i_1}\wedge\cdots\wedge dx_{i_k}. \]
Equation (2.11) implies that $\alpha_{\updot}$ is a map of complexes.
The Dwork trace formula, as formulated by Robba\cite{RO}, then gives
\begin{equation}
L({\bf A}^n/{\bf F}_q,f;t)=\prod_{k=0}^n \det(I-t\alpha_k\mid
\Omega_{C(b)}^k)^{(-1)^{k+1}}.
\end{equation}
From results of Serre\cite{SE} we then get
\begin{equation}
L({\bf A}^n/{\bf F}_q,f;t)=\prod_{k=0}^n \det(I-t\alpha_k\mid
H^k(\Omega_{C(b)}^{\updot},D))^{(-1)^{k+1}},
\end{equation}
where we denote the induced map on cohomology by $\alpha_k$ also.
\section{Filtration by $p$-divisibility}

The $p$-adic Banach space $C(b)$ has a decreasing filtration
$\{F^rC(b)\}_{r=-\infty}^{\infty}$ defined by setting
\[ F^rC(b)=\{\sum_{u\in{\bf N}^n} A_u\pi^{Mb|u|}x^u\in
C(b) \mid A_u\in\pi^r{\mathcal O}_{\tilde{\Omega}_0} \text{ for all
  $u$}\}, \]
where ${\mathcal O}_{\tilde{\Omega}_0}$ denotes the ring of integers
of $\tilde{\Omega}_0$.  We extend this to a filtration on
$\Omega^{\updot}_{C(b)}$ by defining
\[ F^r\Omega^k_{C(b)}=\bigoplus_{1\leq i_1<\cdots<i_k\leq n} F^rC(b)\,
  dx_{i_1}\wedge\cdots\wedge dx_{i_k}. \]
This filtration is exhaustive and separated, i.~e.,
\[ \bigcup_{r\in{\bf Z}}
F^r\Omega^{\updot}_{C(b)}=\Omega^{\updot}_{C(b)} \quad\text{and}\quad
\bigcap_{r\in{\bf Z}}F^r\Omega^{\updot}_{C(b)}=(0). \]
We normalize the $D_i$ so that they respect this filtration.  Put
\[ \epsilon=Mb(d-1)-Md/(p-1), \]
a nonnegative integer.  Then
\[ \pi^{\epsilon}D_i(F^rC(b))\subseteq F^rC(b) \]
and the complexes $(\Omega^{\updot}_{C(b)},D)$,
$(\Omega^{\updot}_{C(b)},\pi^{\epsilon}D)$ have the same cohomology.  

Since $(\Omega^{\updot}_{C(b)},\pi^{\epsilon}D)$ is a filtered
complex, there is an associated spectral sequence.  Its $E_1$-term is
given by
\[ E_1^{r,s}=H^{r+s}(F^r\Omega^{\updot}_{C(b)}/
F^{r+1}\Omega^{\updot}_{C(b)}).  \]
Consider the map $F^0C(b)\rightarrow {\bf F}_q[x_1,\ldots,x_n]$
defined by
\[ \sum_u A_u\pi^{Mb|u|}x^u \mapsto \sum_u
\bar{A}_ux^u, \]
where $\bar{A}_u$ denotes the reduction of $A_u$ modulo the maximal
ideal of ${\mathcal O}_{\tilde{\Omega}_0}$. (Since $A_u\rightarrow 0$
as $u\rightarrow\infty$, the sum on the right-hand side is finite.)
This map induces an isomorphism
\begin{equation}
F^0\Omega^k_{C(b)}/F^1\Omega^k_{C(b)}\simeq \Omega^k_{{\bf
    F}_q[x_1,\ldots,x_n]/{\bf F}_q}.  
\end{equation}
In particular,
\begin{equation}
F^0C(b)/F^1C(b)\simeq{\bf F}_q[x_1,\ldots,x_n].
\end{equation}
We have clearly
\[ \frac{\partial}{\partial x_i}(F^rC(b))\subseteq F^{r+1}C(b), \]
and a calculation show that
\begin{align*} \pi^{\epsilon}H_i&\equiv \pi^{Mb(d-1)}\hat{f}_i \pmod
  {F^1C(b)} \\
 &\equiv \pi^{Mb(d-1)}\hat{f}_i^{(d)} \pmod {F^1C(b)},
\end{align*}
hence under the isomorphism (3.2), the map
\[ \pi^{\epsilon}D_i:F^0C(b)\rightarrow F^0C(b) \]
induces the map ``multiplication by $\partial f^{(d)}/\partial x_i$'' on
${\bf F}_q[x_1,\ldots,x_n]$.  More generally, one sees that under the
isomorphism (3.1), the map
\[ \pi^{\epsilon}D:F^0\Omega^k_{C(b)}\rightarrow
F^0\Omega^{k+1}_{C(b)} \]
induces the map
\[ \phi_{f^{(d)}}:\Omega^k_{{\bf F}_q[x_1,\ldots,x_n]/{\bf
    F}_q}\rightarrow \Omega^{k+1}_{{\bf F}_q[x_1,\ldots,x_n]/{\bf
    F}_q} \]
defined by
\[ \phi_{f^{(d)}}(\omega)=df^{(d)}\wedge\omega, \]
where $df^{(d)}$ denotes the exterior derivative of $f^{(d)}$.  We
have proved that there is an isomorphism of complexes of ${\bf
  F}_q$-vector spaces
\[ (F^0\Omega^{\updot}_{C(b)}/F^1\Omega^{\updot}_{C(b)},\pi^{\epsilon}
D)\simeq (\Omega^{\updot}_{{\bf F}_q[x_1,\ldots,x_n]/{\bf F}_q},
\phi_{f^{(d)}}).  \]
Since multiplication by $\pi^r$ defines an isomorphism of complexes
\[ (F^0\Omega^{\updot}_{C(b)},\pi^{\epsilon}D)\simeq
(F^r\Omega^{\updot}_{C(b)},\pi^{\epsilon}D), \]
we have in fact isomorphisms for all $r\in{\bf Z}$
\begin{equation}
(F^r\Omega^{\updot}_{C(b)}/F^{r+1}\Omega^{\updot}_{C(b)},\pi^{\epsilon}
D)\simeq (\Omega^{\updot}_{{\bf F}_q[x_1,\ldots,x_n]/{\bf F}_q},
\phi_{f^{(d)}}).  
\end{equation}

The complex $(\Omega^{\updot}_{{\bf F}_q[x_1,\ldots,x_n]/{\bf F}_q},
\phi_{f^{(d)}})$ is isomorphic to the Koszul complex on ${\bf
  F}_q[x_1,\ldots,x_n]$ defined by $\{\partial f^{(d)}/\partial
x_i\}_{i=1}^n$.  If we assume
$\{\partial f^{(d)}/\partial x_i\}_{i=1}^n$ form a regular sequence in
${\bf F}_q[x_1,\ldots,x_n]$, we get
\begin{align}
H^i(\Omega^{\updot}_{{\bf F}_q[x_1,\ldots,x_n]/{\bf F}_q},
\phi_{f^{(d)}})& = 0 \qquad \text{for $i\neq n$,} \\
\dim_{{\bf F}_q} H^n(\Omega^{\updot}_{{\bf F}_q[x_1,\ldots,x_n]/{\bf F}_q},
\phi_{f^{(d)}})& =(d-1)^n.
\end{align}
It follows from these equations that
\begin{align}
E_1^{r,s}& =0 \qquad \text{if $r+s\neq n$} \\
\dim_{{\bf F}_q} E_1^{r,s}& =(d-1)^n \qquad \text{if $r+s=n$.}
\end{align}
The first of these equalities implies that all the coboundary maps
$d_1^{r,s}$ are zero, hence the spectral sequence converges weakly, i.~e.,
\[ E_1^{r,s}\simeq F^rH^{r+s}(\Omega^{\updot}_{C(b)},\pi^{\epsilon}D)/ 
F^{r+1}H^{r+s}(\Omega^{\updot}_{C(b)},\pi^{\epsilon}D). \]
This spectral sequence actually converges.  First
observe the following.  Let $x^{\mu_i}$, $i=1,\ldots,(d-1)^n$, be monomials
in $x_1,\ldots,x_n$ such that the cohomology classes $\{[x^{\mu_i}\,
dx_1\wedge\cdots\wedge dx_n]\}_{i=1}^{(d-1)^n}$ form a basis for
$H^n(\Omega^{\updot}_{{\bf F}_q[x_1,\ldots,x_n]/{\bf
    F}_q},\phi_{f^{(d)}})$.  Then the images of the cohomology classes
$\{[\pi^rx^{\mu_i}\,dx_1\wedge\cdots\wedge dx_n]\}_{i=1}^{(d-1)^n}$ in
$E_1^{r,s}$ form a basis for $E_1^{r,s}$ when $r+s=n$.  

\begin{theorem}
  Suppose $\{\partial f^{(d)}/\partial x_i\}_{i=1}^n$ form a regular
  sequence in ${\bf F}_q[x_1,\ldots,x_n]$.  Then \\
  $(1)$ $H^i(\Omega^{\updot}_{C(b)},\pi^{\epsilon}D)=0$ if $i\neq n$, \\
  $(2)$ The cohomology classes $[x^{\mu_i}\, dx_1\wedge\cdots\wedge
  dx_n]$, $i=1,\ldots,(d-1)^n$, form a basis for
  $H^n(\Omega^{\updot}_{C(b)},\pi^{\epsilon}D)$.
\end{theorem}

{\it Proof}.  Suppose $i\neq n$ and let $\eta\in\Omega^i_{C(b)}$ with
$\pi^\epsilon D(\eta)=0$.  For some $r$ we have $\eta\in
F^r\Omega^i_{C(b)}$.  Equations (3.3) and (3.4) then imply that
\[ \eta=\pi\eta_1+\pi^{\epsilon}D(\zeta_1) \]
with $\eta_1\in F^r\Omega^i_{C(b)}$ and $\zeta_1\in
F^r\Omega^{i-1}_{C(b)}$.  Suppose that for some $t\geq 1$ we have
found $\eta_t\in F^r\Omega^i_{C(b)}$ and $\zeta_t\in
F^r\Omega^{i-1}_{C(b)}$ such that 
\begin{equation}
\eta=\pi^t\eta_t+\pi^{\epsilon}D(\zeta_t)
\end{equation}
and such that
\[ \zeta_t-\zeta_{t-1}\in F^{r+t-1}\Omega^{i-1}_{C(b)}. \]
Applying $\pi^{\epsilon}D$ to both sides of (3.9) gives
\[ \pi^{t+\epsilon}D(\eta_t)=0, \]
hence $\pi^{\epsilon}D(\eta_t)=0$ since multiplication by $\pi$ is
injective on $\Omega^i_{C(b)}$.  Equations (3.3) and (3.4) give
\[ \eta_t=\pi\eta_{t+1}+\pi^{\epsilon}D(\zeta'_{t+1}), \]
with $\eta_{t+1}\in F^r\Omega^i_{C(b)}$ and $\zeta'_{t+1}\in
F^r\Omega^{i-1}_{C(b)}$.  If we put
$\zeta_{t+1}=\zeta_t+\pi^t\zeta'_{t+1}$, then substitution into (3.9)
gives 
\[ \eta=\pi^{t+1}\eta_{t+1}+\pi^{\epsilon}D(\zeta_{t+1}) \]
with
\[ \zeta_{t+1}-\zeta_t\in F^{r+t}\Omega^{i-1}_{C(b)}. \]
It is now clear that the sequence $\{\zeta_t\}_{t=1}^{\infty}$
converges to an element $\zeta\in F^r\Omega^{i-1}_{C(b)}$ such that
$\eta=\pi^{\epsilon}D(\zeta)$.  This proves the first assertion.

It follows easily from (3.3) that the $\{[x^{\mu_i}\,dx_1\wedge\cdots\wedge
dx_n]\}_{i=1}^{(d-1)^n}$ are linearly independent, hence it suffices to
show that they span $H^n(\Omega^{\updot}_{C(b)},\pi^{\epsilon}D)$.
Let $\eta\in F^r\Omega^n_{C(b)}$.  From (3.3) we have
\[ \eta=\sum_{i=1}^{(d-1)^n} c_i^{(1)}x^{\mu_i}\,dx_1\wedge\cdots\wedge
dx_n+\pi\eta_1+\pi^{\epsilon}D(\zeta_1), \]
where $c_i^{(1)}\in\tilde{\Omega}_0$, $c_i^{(1)}x^{\mu_i}\in F^rC(b)$,
$\eta_1\in F^r\Omega^n_{C(b)}$, $\zeta_1\in F^r\Omega^{n-1}_{C(b)}$.
Suppose we can write
\begin{equation}
\eta=\sum_{i=1}^{(d-1)^n} c_i^{(t)}x^{\mu_i}\,dx_1\wedge\cdots\wedge
dx_n+\pi^t\eta_t+\pi^{\epsilon}D(\zeta_t)
\end{equation}
with $c_i^{(t)}\in\tilde{\Omega}_0$, $c_i^{(t)}x^{\mu_i}\in F^rC(b)$,
$\eta_t\in F^r\Omega^n_{C(b)}$, and $\zeta_t\in
F^r\Omega^{n-1}_{C(b)}$ such that 
\begin{align*}
(c_i^{(t)}-c_i^{(t-1)})x^{\mu_i}&\in F^{r+t-1}C(b) \\
\zeta_t-\zeta_{t-1}&\in F^{r+t-1}\Omega^{n-1}_{C(b)}.
\end{align*}
By (3.3) we have
\[ \eta_t=\sum_{i=1}^{(d-1)^n} c'_ix^{\mu_i}\,dx_1\wedge\cdots\wedge
dx_n+\pi\eta_{t+1}+\pi^{\epsilon}D(\zeta'_t), \]
where $c'_i\in\tilde{\Omega}_0$, $c'_ix^{\mu_i}\in F^rC(b)$,
$\eta_{t+1}\in F^r\Omega^n_{C(b)}$, $\zeta'_t\in
F^r\Omega^{n-1}_{C(b)}$.  If we put $c_i^{(t+1)}=c_i^{(t)}+\pi^tc'_i$
and $\zeta_{t+1}=\zeta_t+\pi^t\zeta'_t$, then
\[ \eta=\sum_{i=1}^{(d-1)^n} c_i^{(t+1)}x^{\mu_i}\,dx_1\wedge\cdots\wedge
dx_n+\pi^{t+1}\eta_{t+1}+\pi^{\epsilon}D(\zeta_{t+1}) \]
with
\begin{align*}
(c_i^{(t+1)}-c_i^{(t)})x^{\mu_i}&\in F^{r+t}C(b) \\
\zeta_{t+1}-\zeta_t&\in F^{r+t}\Omega^{n-1}_{C(b)}.
\end{align*}
It follows that the sequences $\{c_i^{(t)}\}_{t=1}^{\infty}$ and
$\{\zeta_t\}_{t=1}^{\infty}$ converge, say, $c_i^{(t)}\rightarrow
c_i\in\tilde{\Omega}_0$,$\zeta_t\rightarrow\zeta\in
F^r\Omega^{n-1}_{C(b)}$, and that these limits satisfy
\[ \eta=\sum_{i=1}^{(d-1)^n} c_ix^{\mu_i}\,dx_1\wedge\cdots\wedge
dx_n+\pi^{\epsilon}D(\zeta) \]
with $c_ix^{\mu_i}\in F^rC(b)$.  This completes the proof of the
second assertion. 

The following result is a consequence of the proof of Theorem 3.8.  
\begin{proposition}
Under the hypothesis of Theorem $3.8$, if $\eta\in
F^r\Omega^n_{C(b)}$, then there exist
$\{c_i\}_{i=1}^{(d-1)^n}\subseteq \tilde{\Omega}_0$ such that in
$H^n(\Omega^{\updot}_{C(b)},\pi^{\epsilon}D)$ we have
\[ [\eta]=\sum_{i=1}^{(d-1)^n}[c_ix^{\mu_i}\,dx_1\wedge\cdots\wedge dx_n], \]
where $c_ix^{\mu_i}\in F^rC(b)$ for $i=1,\ldots,(d-1)^n$.
\end{proposition}

\section{$p$-adic estimates}

It follows from (2.14) and Theorem 3.8 that
\begin{equation}
L({\bf A}^n,f;t)^{(-1)^{n+1}}=\det(I-t\alpha_n\mid
H^n(\Omega^{\updot}_{C(b)},D)) 
\end{equation}
is a polynomial of degree $(d-1)^n$ (by \cite{RO}, zero is not an
eigenvalue of $\alpha_n$).  We estimate its $p$-adic Newton polygon.
Note that
\[ H^n(\Omega^{\updot}_{{\bf F}_q[x_1,\ldots,x_n]/{\bf
    F}_q},\phi_{f^{(d)}}) \simeq {\bf F}_q[x_1,\ldots,x_n]/(\partial
f^{(d)}/\partial x_1,\ldots,\partial f^{(d)}/\partial x_n) \]
is a graded ${\bf F}_q[x_1,\ldots,x_n]$-module.  Let
$H^n(\Omega^{\updot}_{{\bf F}_q[x_1,\ldots,x_n]/{\bf
    F}_q},\phi_{f^{(d)}})^{(m)}$ denote its homogeneous component of
degree $m$.  It follows from (3.4) that its Hilbert-Poincare series is
$(1+t+\cdots+t^{d-2})^n$.  Write
\begin{equation}
(1+t+\cdots+t^{d-2})^n =\sum_{m=0}^{n(d-2)} U_mt^m,
\end{equation}
so that
\[ U_m=\dim_{{\bf F}_q} H^n(\Omega^{\updot}_{{\bf F}_q[x_1,\ldots,x_n]/{\bf
    F}_q},\phi_{f^{(d)}})^{(m)}. \]
Equivalently, 
\[ U_m={\rm card} \{x^{\mu_i}\mid |\mu_i|=m\}. \]
\begin{theorem}
  Suppose $\{\partial f^{(d)}/\partial x_i\}_{i=1}^n$ form a regular
    sequence in ${\bf F}_q[x_1,\ldots,x_n]$.  Then the Newton polygon
    of $L({\bf A}^n,f;t)^{(-1)^{n+1}}$ with respect to the valuation
    ``${\rm ord}_q$'' lies on or above the Newton polygon with respect
    to the valuation ``${\rm ord}_q$'' of the polynomial
\[ \prod_{m=0}^{n(d-2)} (1-q^{(m+n)/d}t)^{U_m}. \]
\end{theorem}

We begin with a reduction step.  Let $\beta_n$ be the endomorphism of
$H^n(\Omega^{\updot}_{C(b)},D)$ constructed from $\beta$ as $\alpha_n$
was constructed from $\alpha$, i.~e.,
\[ \beta_n(\xi\,dx_1\wedge\cdots\wedge dx_n)=\frac{1}{x_1\cdots
  x_n}\beta(x_1\cdots x_n\xi)\,dx_1\wedge\cdots\wedge dx_n. \]
Then $\beta_n$ is an $\tilde{\Omega}_1$-linear endomorphism of
$H^n(\Omega^{\updot}_{C(b)},D)$ and $\alpha_n=(\beta_n)^a$.  By
\cite[Lemma 7.1]{DW2}, we have 
the following.
\begin{lemma}  
  The Newton polygon of $\det_{\tilde{\Omega}_0}(I-t\alpha_n\mid
  H^n(\Omega^{\updot}_{C(b)},D))$ with respect to the valuation
  ``${\rm ord}_q$'' is obtained from the Newton polygon of
  $\det_{\tilde{\Omega}_1}(I-t\beta_n\mid
  H^n(\Omega^{\updot}_{C(b)},D))$ with respect to the valuation
  ``ord'' by shrinking ordinates and abscissas by a factor of $1/a$.  
\end{lemma}

Let $\{\gamma_j\}_{j=1}^a$ be an integral basis for $\tilde{\Omega}_0$
over $\tilde{\Omega}_1$.  Then under the hypothesis of Theorem 3.8,
the cohomology classes
\[ [\gamma_jx^{\mu_i}\,dx_1\wedge\cdots\wedge dx_n],\qquad
\text{$i=1,\ldots,(d-1)^n$, $j=1,\ldots,a$,} \]
form a basis for $H^n(\Omega^{\updot}_{C(b)},D)$ as
$\tilde{\Omega}_1$-vector space.  We estimate $p$-adically the entries
of the matrix of $\beta_n$ with respect to a certain normalization of
this basis, namely, we set
\[ \xi(i,j)=(\pi^{Mb/p})^{|\mu_i|+n}\gamma_jx^{\mu_i} \]
and use the cohomology classes $[\xi(i,j)\,dx_1\wedge\cdots\wedge
dx_n]$.  This normalization is chosen so that
\[ x_1\cdots x_n\xi(i,j)\in C(b/p,0), \]
hence
\[ \beta(x_1\cdots x_n\xi(i,j))\in C(b,0) \]
and
\[ \frac{1}{x_1\cdots x_n}\beta(x_1\cdots x_n\xi(i,j))\in
\pi^{Mbn}C(b,0). \]
This says that 
\[ \beta_n(\xi(i,j)\,dx_1\wedge\cdots\wedge dx_n)\in
F^{Mbn}\Omega^n_{C(b)}, \]
hence by Proposition 3.11 and the properties of an integral basis we
have
\[ [\beta_n(\xi(i,j)\,dx_1\wedge\cdots\wedge dx_n)]=\sum_{i',j'}
A(i',j';i,j)[\gamma_{j'}x^{\mu_{i'}}\,dx_1\wedge\cdots\wedge dx_n] \]
with $A(i',j';i,j)\gamma_{j'}x^{\mu_{i'}}\in F^{Mbn}C(b)$, i.~e.,
$A(i',j';i,j)\in \pi^{Mb(|\mu_i|+n)}{\mathcal O}_{\tilde{\Omega}_0}$.
This may be rewritten as
\[ [\beta_n(\xi(i,j)\,dx_1\wedge\cdots\wedge dx_n)]=\sum_{i',j'}
B(i',j';i,j)[\xi(i',j')\,dx_1\wedge\cdots\wedge dx_n] \]
with
\[ B(i',j';i,j)\in\pi^{Mb(|\mu_{i'}|+n)(1-1/p)}{\mathcal
    O}_{\tilde{\Omega}_0}, \]
i.~e., the $(i',j')$-row of the matrix $B(i',j';i,j)$ of $\beta_n$
with respect to the basis $\{[\xi(i,j)\,dx_1\wedge\cdots\wedge
dx_n)]\}_{i,j}$ is divisible by
\[ \pi^{Mb(|\mu_{i'}|+n)(1-1/p)}. \]
This implies that $\det_{\tilde{\Omega}_1}(I-t\beta_n\mid
H^n(\Omega^{\updot}_{C(b)},D))$ has Newton polygon (with respect to
the valuation ``ord'') lying on or above the Newton polygon (with
respect to the valuation ``ord'') of the polynomial
\[ \prod_{m=0}^{n(d-2)} (1-\pi^{Mb(m+n)(1-1/p)}t)^{aU_m}. \]
But $\det_{\tilde{\Omega}_1}(I-t\beta_n\mid
H^n(\Omega^{\updot}_{C(b)},D))$ is independent of $b$, so we may take
the limit as $b\rightarrow p/(p-1)$ to conclude that its Newton
polygon lies on or above the Newton polygon of
\[ \prod_{m=0}^{n(d-2)} (1-p^{(m+n)/d}t)^{aU_m}. \]
Theorem 4.3 now follows from Lemma 4.4.

Let $\{\rho_i\}_{i=1}^{(d-1)^n}$ be the reciprocal roots of $L({\bf
  A}^n,f;t)^{(-1)^{n+1}}$ and put
\[ \Lambda(f)=\prod_{i=1}^{(d-1)^n} \rho_i\in{\bf Q}(\zeta_p). \]
Theorem 4.3 implies that
\[ {\rm ord}_q\;\Lambda(f)\geq \frac{1}{d}\sum_{m=0}^{n(d-2)}(m+n)U_m. \]
But it follows from (4.2) evaluated at $t=1$ that
\[ \sum_{m=0}^{n(d-2)} U_m=(d-1)^n \]
and from the derivative of (4.2) evaluated at $t=1$ that
\[ \sum_{m=0}^{n(d-2)} mU_m=n(d-1)^n(d-2)/2. \]
We thus get the following.
\begin{corollary}
Under the hypothesis of Theorem $4.3$,
\[ {\rm ord}_q\;\Lambda(f)\geq \frac{n(d-1)^n}{2}. \]
\end{corollary}

It can be proved directly by $p$-adic methods that equality holds in
Corollary 4.5.  We shall derive this equality in the next section by
$l$-adic methods.

\section{$l$-adic cohomology}

Let $l$ be a prime, $l\neq p$.  There exists a lisse, rank-one,
$l$-adic \'{e}tale sheaf ${\mathcal L}_{\Psi}(f)$ on ${\bf A}^n$ with
the property that
\begin{equation}
L({\bf A}^n,f;t)=L({\bf A}^n,{\mathcal L}_{\Psi}(f);t),
\end{equation}
where the right-hand side is a Grothendieck $L$-function.  By
Grothendieck's Lefschetz trace formula,
\begin{equation}
L({\bf A}^n,f;t)=\prod_{i=0}^{2n} \det(I-tF\mid H^i_c({\bf
  A}^n\times_{{\bf F}_q}\bar{{\bf F}}_q,{\mathcal
  L}_{\Psi}(f)))^{(-1)^{i+1}},
\end{equation}
where $H_c^i$ denotes $l$-adic cohomology with proper supports and $F$
is the Frobenius endomorphism.  The $H^i_c({\bf A}^n\times_{{\bf
    F}_q}\bar{{\bf F}}_q,{\mathcal L}_{\Psi}(f))$ are
finite-dimensional vector spaces over a finite extension $K_l$ of
${\bf Q}_l$ containing the $p$-th roots of unity.  We combine Theorem
3.8 and Corollary 4.5 with general results of Deligne\cite{DE2} to
prove the following theorem of Deligne\cite[Th\'{e}or\`{e}me
8.4]{DE1}.  
\begin{theorem}
Suppose $(p,d)=1$ and $f^{(d)}=0$ defines a smooth hypersurface
in~${\bf P}^{n-1}$.  Then \\ 
$(1)$ $H^i_c({\bf A}^n\times_{{\bf F}_q}\bar{{\bf F}}_q,{\mathcal
  L}_{\Psi}(f))=0$ if $i\neq n$, \\
$(2)$ $\dim_{K_l} H^n_c({\bf A}^n\times_{{\bf F}_q}\bar{{\bf
    F}}_q,{\mathcal L}_{\Psi}(f))=(d-1)^n, \\
$(3) $H^n_c({\bf A}^n\times_{{\bf F}_q}\bar{\bf F}_q,{\mathcal
  L}_{\Psi}(f))$ is pure of weight $n$.
\end{theorem}

{\it Proof}.  We consider the theorem to be known for $n=1$ and prove
it for general $n\geq 2$ by induction.  For $\lambda\in\bar{\bf F}_q$,
set
\[ f_{\lambda}(x_1,\ldots,x_{n-1})=f(x_1,\ldots,x_{n-1},\lambda)\in
{\bf F}_q(\lambda)[x_1,\ldots,x_{n-1}]. \]
Since the generic hyperplane section of a smooth variety is smooth, we
may assume, after a coordinate change if necessary, that the
hyperplane $x_n=0$ intersects the variety $f^{(d)}=0$ transversally in
${\bf P}^{n-1}$.  Thus $f^{(d)}(x_1,\ldots,x_{n-1},0)=0$ defines a
smooth hypersurface in ${\bf P}^{n-2}$.  But
\[ f_{\lambda}^{(d)}=f^{(d)}(x_1,\ldots,x_{n-1},0), \]
so by the induction hypothesis the conclusions of the theorem are true
for all $f_{\lambda}$.

Consider the morphism of ${\bf F}_q$-schemes $\sigma:{\bf
  A}^n\rightarrow{\bf A}^1$ which is projection onto the $n$-th
  coordinate.  The Leray spectral sequence for the composition of
  $\sigma$ with the structural morphism ${\bf A}^1\rightarrow{\rm
  Spec}({\bf F}_q)$ is
\begin{equation}
H^i_c({\bf A}^1\times_{{\bf F}_q}\bar{\bf F}_q,R^j\sigma_!({\mathcal
  L}_{\Psi}(f)))\Rightarrow H^{i+j}_c({\bf A}^n\times_{{\bf
  F}_q}\bar{{\bf F}}_q,{\mathcal L}_{\Psi}(f)).
\end{equation}
Proper base change implies that for $\lambda\in\bar{\bf F}_q$ and
$\bar{\lambda}$ a geometric point over $\lambda$
\begin{equation}
(R^j\sigma_!({\mathcal L}_{\Psi}(f)))_{\bar{\lambda}}=H^j_c(
  \sigma^{-1}(\lambda)\times_{{\bf F}_q(\lambda)}\bar{{\bf
  F}}_q,{\mathcal L}_{\Psi}(f_{\lambda})).
\end{equation}
Applying the induction hypothesis to $f_{\lambda}$ shows that the
right-hand side of (5.5) vanishes for all $\lambda\in\bar{\bf F}_q$ if
$j\neq n-1$.  It follows that the Leray spectral sequence collapses
and we get
\begin{equation}
H^i_c({\bf A}^1\times_{{\bf F}_q}\bar{\bf F}_q,R^{n-1}\sigma_!({\mathcal
  L}_{\Psi}(f)))=H^{i+n-1}_c({\bf A}^n\times_{{\bf F}_q}\bar{{\bf
  F}}_q,{\mathcal L}_{\Psi}(f)).
\end{equation}
Since $\dim{\bf A}^1=1$, the left-hand side of (5.6) can be nonzero only
for $i=0,1,2$.  However, $H^{n-1}_c({\bf A}^n\times_{{\bf F}_q}\bar{{\bf
  F}}_q,{\mathcal L}_{\Psi}(f))=0$ because ${\bf A}^n$ is smooth,
affine, of dimension~$n$, and ${\mathcal L}_{\Psi}(f)$ is lisse on
${\bf A}^n$.  This proves that $H^i_c({\bf A}^n\times_{{\bf F}_q}\bar{{\bf
  F}}_q,{\mathcal L}_{\Psi}(f))=0$ except possibly for $i=n,n+1$.

By (5.2) we then have
\begin{equation}
L({\bf A}^n,f;t)^{(-1)^{n+1}}=\frac{\det(I-tF\mid H^n_c({\bf
  A}^n\times_{{\bf F}_q}\bar{{\bf F}}_q,{\mathcal L}_{\Psi}(f)))}
{\det(I-tF\mid H^{n+1}_c({\bf A}^n\times_{{\bf F}_q}\bar{{\bf
  F}}_q,{\mathcal L}_{\Psi}(f)))}. 
\end{equation}
Since ${\mathcal L}_{\Psi}(f)$ is pure of weight 0, Deligne's
fundamental theorem\cite{DE2} tells us that $H^n_c({\bf
  A}^n\times_{{\bf F}_q}\bar{{\bf F}}_q,{\mathcal L}_{\Psi}(f))$ is
mixed of weights $\leq n$.  Equation (5.5) and the induction hypothesis
applied to $f_{\lambda}$ tell us that $R^{n-1}\sigma_!({\mathcal
  L}_{\Psi}(f))$ is pure of weight $n-1$ and that all fibers of
$R^{n-1}\sigma_!({\mathcal L}_{\Psi}(f))$ 
have the same rank, namely, $(d-1)^{n-1}$.  It follows from
Katz\cite[Corollary 6.7.2]{KA} that $R^{n-1}\sigma_!({\mathcal
  L}_{\Psi}(f))$ is lisse on ${\bf A}^1$.  Equation (5.6) with $i=2$ now
implies, by Deligne\cite[Corollaire 1.4.3]{DE2}, that $H^{n+1}_c({\bf
  A}^n\times_{{\bf F}_q}\bar{{\bf F}}_q,{\mathcal L}_{\Psi}(f))$ is
pure of weight $n+1$, hence there can be no cancellation on the
right-hand side of (5.7).  However, Theorem 3.8 implies that $L({\bf
  A}^n,f;t)^{(-1)^{n+1}}$ is a polynomial of degree $(d-1)^n$, so we
must have
\[ H^{n+1}_c({\bf A}^n\times_{{\bf F}_q}\bar{{\bf F}}_q,{\mathcal
  L}_{\Psi}(f))=0 \]
and
\[ \dim_{K_l} H^n_c({\bf A}^n\times_{{\bf F}_q}\bar{{\bf F}}_q,{\mathcal
  L}_{\Psi}(f))=(d-1)^n. \]
This establishes the first two assertions of the theorem.

To prove the last assertion of the theorem, note that $|\rho_i|\leq
q^{n/2}$ for every $i$ and every archimedian absolute value since
$H^n_c({\bf A}^n\times_{{\bf F}_q}\bar{{\bf F}}_q,{\mathcal 
  L}_{\Psi}(f))$ is mixed of weights $\leq n$.  Thus we have
\begin{equation}
|\Lambda(f)|\leq q^{n(d-1)^n/2}
\end{equation}
for every archimedian absolute value on ${\bf Q}(\zeta_p)$.  By
Corollary 4.5, we have
\begin{equation}
|\Lambda(f)|_p\leq q^{-n(d-1)^n/2}
\end{equation}
for every normalized archimedian absolute value on ${\bf Q}(\zeta_p)$
lying over $p$, and it is well-known that $|\rho_i|_{p'}=1$ for every
nonarchimedian absolute value lying over any prime $p'\neq p$.  It
then follows from the product formula for ${\bf Q}(\zeta_p)$ that
equality holds in (5.8) (and also in (5.9)), which implies the last
assertion of the theorem.

\section{Proof of Theorem 1.4}

It remains to consider the case where $p$ divides $d$.  The Euler
relation becomes
\[ \sum_{i=1}^n x_i\frac{\partial f^{(d)}}{\partial x_i}=0. \]
The regular sequence hypothesis then implies that
\[ x_i\in \biggl(\frac{\partial f^{(d)}}{\partial
  x_1},\ldots,\widehat{\frac{\partial f^{(d)}}{\partial x_i}},\ldots,
  \frac{\partial f^{(d)}}{\partial x_n}\biggr), \]
hence there is an equality of ideals of ${\bf F}_q[x_1,\ldots,x_n]$
\begin{equation}
(x_1,\ldots,x_n)=\biggl(\frac{\partial f^{(d)}}{\partial
x_1},\ldots,\frac{\partial f^{(d)}}{\partial x_n}\biggr). 
\end{equation}
Conversely, if (6.1) holds, then $\{\partial f^{(d)}/\partial
x_i\}_{i=1}^n$ is a regular sequence.  Equation (6.1) implies that
$d=2$, hence $p=2$ as well, thus $f$ is a quadratic polynomial in
characteristic $2$.  We may assume $f$ contains no terms of the form
$x_i^2$ by the following elementary lemma.

Let $\zeta_p$ be a primitive $p$-th root of unity.  Since $\Psi$ is a
nontrivial additive character of ${\bf F}_q$, there exists a nonzero
$b\in{\bf F}_q$ such that
\begin{equation}
\Psi(x)=\zeta_p^{{\rm Tr}_{{\bf F}_q/{\bf F}_p}(bx)}.
\end{equation}
\begin{lemma}
Let $a\in{\bf F}_q$, $a\neq 0$, and choose $c\in{\bf F}_q$ such that
$c^p=(ab)^{-1}$.  Then
\[ \sum_{x_1,\ldots,x_n\in{\bf F}_q} \Psi(f(x_1,\ldots,x_n)+ax_n^p)=
\sum_{x_1,\ldots,x_n\in{\bf F}_q}
\Psi(f(x_1,\ldots,x_n)+ac^{p-1}x_n). \]
\end{lemma}

{\it Proof}.  Making the change of variable $x_n\mapsto cx_n$, the sum
becomes 
\begin{multline*}
\sum_{x_1,\ldots,x_n\in{\bf F}_q}
\Psi(f(x_1,\ldots,x_{n-1},cx_n)+b^{-1}x_n^p) = \\
\sum_{x_1,\ldots,x_n\in{\bf F}_q}
\Psi(f(x_1,\ldots,x_{n-1},cx_n)+b^{-1}x_n) \Psi(b^{-1}(x_n^p-x_n)). 
\end{multline*}
But by (6.2), 
\begin{align*}
\Psi(b^{-1}(x_n^p-x_n))& = \zeta_p^{{\rm Tr}_{{\bf F}_q/{\bf
      F}_p}(x_n^p-x_n)} \\
 & = 1
\end{align*}
since ${\rm Tr}_{{\bf F}_q/{\bf F}_p}(x_n^p-x_n)=0$ for all
$x_n\in{\bf F}_q$.  Making the change of variable $x_n\mapsto
c^{-1}x_n$ now gives the lemma.

By the lemma, we may assume our quadratic polynomial $f$ has the form
\[ f=\sum_{1\leq i<j\leq n}a_{ij}x_ix_j + \sum_{k=1}^n b_kx_k +c, \]
where $a_{ij},b_k,c\in{\bf F}_q$.  This gives
\[ f^{(2)}=\sum_{1\leq i<j\leq n}a_{ij}x_ix_j. \]
Let $A=(A_{ij})$ be the $n\times n$ matrix defined by
\[ A_{ij}=\begin{cases}
a_{ij}& \text{if $i<j$} \\
0& \text{if $i=j$} \\
a_{ji}& \text{if $i>j$.}
\end{cases} \]
Thus $A$ is a symmetric matrix with zeros on the diagonal.  One checks
that
\[ \begin{bmatrix} \frac{\partial f^{(2)}}{\partial x_1} \\
  \hdotsfor{1} \\ \frac{\partial f^{(2)}}{\partial x_n} \end{bmatrix} =
  A \begin{bmatrix} x_1 \\ \hdotsfor{1} \\ x_n \end{bmatrix}, \]
therefore (6.1) holds if and only if $\det A\neq 0$.  

We now evaluate the exponential sum
\begin{equation}
\sum_{x_1,\ldots,x_n\in{\bf F}_q} \Psi\biggl(\sum_{1\leq i<j\leq
  n}a_{ij}x_ix_j + \sum_{k=1}^n b_kx_k +c\biggr). 
\end{equation}
\begin{proposition}
If $n$ is odd, then $(6.1)$ cannot hold.  If $n$ is even and $(6.1)$
holds, then the sum $(6.4)$ equals $\zeta q^{n/2}$, where $\zeta$ is a
root of unity.
\end{proposition}

{\it Proof}.  If $n=1$, then $\det A=0$, so (6.1) cannot hold.  If
$n=2$, then $\det A\neq 0$ if and only if $a_{12}\neq 0$.  It is then
easy to check that the sum (6.4) equals 
\[ \Psi\biggl(\frac{b_1b_2}{a_{12}}+c\biggr)q. \]
Thus the proposition holds for $n=1,2$.  Suppose $n\geq 3$.  The sum
(6.4) can be rewritten as
\[ \sum_{x_1,\ldots,x_{n-1}\in{\bf F}_q} \Psi\biggl(\sum_{1\leq i<j\leq
  n-1}a_{ij}x_ix_j + \sum_{k=1}^{n-1} b_kx_k +c\biggr)\sum_{x_n\in{\bf F}_q}
  \Psi\biggl(\biggl(\sum_{i=1}^{n-1} a_{in}x_i+b_n\biggr)x_n\biggr). \]
But
\[ \sum_{x_n\in{\bf F}_q}\Psi\biggl(\biggl(\sum_{i=1}^{n-1}
  a_{in}x_i+b_n\biggr)x_n\biggr) =
  \begin{cases} 0& \text{if $\sum_{i=1}^{n-1} a_{in}x_i+b_n\neq 0$} \\
  q & \text{if $\sum_{i=1}^{n-1} a_{in}x_i+b_n=0$}, \end{cases} \]
hence (6.4) equals
\begin{equation}
q\sum_{\substack{x_1,\ldots,x_{n-1}\in{\bf F}_q\\
    \sum_{i=1}^{n-1} a_{in}x_i+b_n=0}} \Psi\biggl(\sum_{1\leq i<j\leq
    n-1}a_{ij}x_ix_j + \sum_{k=1}^{n-1} b_kx_k +c\biggr). 
\end{equation}
Since we are assuming $A$ is invertible, some $a_{in}$ must be
nonzero, say, $a_{n-1,n}\neq 0$.  By making the change of variable
$x_{n-1}\mapsto (a_{n-1,n})^{-1}x_{n-1}$, we may assume $a_{n-1,n}=1$.
Solving $a_{1n}x_1+\cdots+a_{n-1,n}x_{n-1}+b_n=0$ for $x_{n-1}$ and
substituting into the expression in the additive character, we see
that (6.6) equals
\begin{equation}
q\sum_{x_1,\ldots,x_{n-2}\in{\bf F}_q} \Psi\biggl(\sum_{1\leq i<j\leq
    n-2}a'_{ij}x_ix_j + \sum_{k=1}^{n-2} b'_kx_k +c\biggr),
\end{equation}
where
\[ a'_{ij}=a_{ij}+a_{i,n-1}a_{jn}+a_{j,n-1}a_{in}. \]
Let $A'=(A'_{ij})$ be the $(n-2)\times(n-2)$ matrix constructed from
the $a'_{ij}$ as $A$ was constructed from the $a_{ij}$.  We explain
the connection between $A$ and $A'$.  Let $\tilde{A}$ be the $n\times
n$ matrix obtained from $A$ by replacing row $i$ by
\[ {\rm row}\;i+a_{in}({\rm row}\;n-1)+a_{i,n-1}({\rm row}\;n) \]
for $i=1,\ldots,n-2$.  Keeping in mind that $a_{n-1,n}=1$, we see that
\[ \tilde{A}=\left[ \begin{array}{ccc|cc} 
 & & & 0 & 0 \\ & A' & & \vdots & \vdots \\ & & & 0 & 0 \\ \hline
a_{1,n-1} & \cdots & a_{n-2,n-1} & 0 & 1 \\
a_{1n} & \cdots & a_{n-2,n} & 1 & 0 \end{array} \right]. \]
In particular, $\det A'=\det A$.  We can repeat this procedure
starting with the sum~(6.7) and continue until we are reduced to the
one or two variable case, according to whether $n$ is odd or even.  If
$n$ is odd, this implies $\det A=0$, a contradiction.  Thus there does
not exist a quadratic polynomial $f$ satisfying (6.1) in this case.
If $n$ is even, this shows that (6.4) equals $q^{n/2}$ times a root of
unity, which is the desired result.

A straightforward calculation using Proposition 6.5 then shows that
the corresponding $L$-function has the form asserted in Theorem 1.4.

\end{document}